\documentclass[12pt,a4paper,twoside]{amsart}
\usepackage{amsmath,amssymb, MnSymbol, amsthm, stmaryrd, MnSymbol, hyperref, amsfonts, bbding, pifont, mathrsfs}

\usepackage{graphicx}
\usepackage{xcolor,color}
\usepackage{enumerate}
\usepackage[utf8]{inputenc}
\usepackage[textsize=small]{todonotes}

\newtheorem{thm}{Theorem}[section]
\newtheorem{cor}[thm]{Corollary}
\newtheorem{prop}[thm]{Proposition}
\newtheorem{lem}[thm]{Lemma}
\newtheorem{fact}[thm]{Fact}

\theoremstyle{definition}
\newtheorem{dfn}[thm]{Definition}

\newcommand\gl{\mathsf{GL}}
\newcommand\lang{\mathcal L}
\newcommand\model{\mathfrak M}
\newcommand\inter[1]{\llbracket #1\rrbracket}
\newcommand\set[1]{\{#1\}}
\newcommand\tup[1]{\langle #1\rangle}
\newcommand\cof{\operatorname{cf}}

\newcommand\mc{\mathfrak c}
\newcommand\pp{\mathfrak p}
\newcommand\Fin{\operatorname{Fin}}
\newcommand\diag{\operatorname{Diag}}
\newcommand\Th{\operatorname{Th}}
\newcommand\restr{\mathbin{\upharpoonright}}
\newcommand\substar{\subseteq^{*}}
\renewcommand\root{\operatorname{root}}

\newcommand{\lbr}{\llbracket}
\newcommand{\rbr}{\rrbracket}

\title[strong completeness below continuum]{Strong Completeness of Provability Logic for Languages of size at most continuum}
\date{}

\author[M. Golshani]{Mohammad Golshani}

\address{Mohammad Golshani, School of Mathematics, Institute for Research in Fundamental Sciences (IPM), P.O.\ Box:
	19395--5746, Tehran, Iran.}

\email{golshani.m@gmail.com}

\author[R. Zoghifard]{Reihane Zoghifard}

\address{Reihane Zoghifard, School of Mathematics, Institute for Research in Fundamental Sciences (IPM), P.O.\ Box:
	19395--5746, Tehran, Iran.}

\email{r.zoghi@gmail.com}

\begin{document}


\begin{abstract}
We investigate the strong completeness of $\mathrm{GL}(\mathbb P)$ for an infinite set of propositional variables $\mathbb P$ over ordinal spaces equipped with the interval topology $\mathcal{I}_1$. For every infinite cardinal $\kappa$, we consider the set
\[
\Gamma_\kappa=
\{\Diamond p_\xi:\xi<\kappa\}
\cup
\{\Box(p_\xi\to\Diamond p_\eta):\xi<\eta<\kappa\},
\]
and prove that $\Gamma_\kappa$ has a model on an ordinal space with $\mathcal{I}_1$ if and only if
$\kappa\nrightarrow_{\mathrm{ipath}}(\omega)^2_{\omega,<\omega}$.
We also prove that if $|\mathbb P|<\mathfrak p$, then $\mathrm{GL}(\mathbb P)$ is strongly complete with respect to $(\omega^\omega+1,\mathcal I_1)$.
Finally, we establish the strong completeness of $\mathrm{GL}(\mathbb P)$ with respect to the same space whenever $|\mathbb P|=\aleph_1$.
\end{abstract}

\maketitle

\section{Introduction}

This paper continues the study of strong completeness of the G\"odel--L\"ob provability logic for ordinal spaces over uncountable languages initiated by \cite{aguilera2017strong, GZS}.

For a given set of propositional variables $\mathbb{P}$,
let $\mathcal{L}(\mathbb{P})$ denote the set of modal formulas in the language containing $\mathbb{P}$, formed using the Boolean connectives and the modal operator $\Box$.
The G\"odel--L\"ob provability logic over $\mathbb{P}$, denoted by $\gl(\mathbb{P})$, is the propositional modal logic
axiomatised by the following axiom schemata
\begin{itemize}
	\item $\text{K.} \Box(\varphi\to \psi) \to (\Box\varphi\to \Box\psi)$,
	\item $\text{L\"ob.} \Box (\Box\varphi\to\varphi) \to\Box\varphi$,
\end{itemize}
with Modus Ponens and Necessitation $\frac{\varphi}{\Box\varphi}$ as inference rules.
When $\mathbb{P}$ is countable, we simply write $\gl$ instead of $\gl(\mathbb{P})$.

Segerberg  \cite{segerberg} studied the Kripke semantics of $\gl$ and proved its completeness with respect to the class of all finite transitive and conversely well-founded trees.
However, $\gl$ is not strongly complete with respect to its Kripke semantics.
Later, Esakia \cite{esakia} observed the analogy between the behaviour of the modal operator $ \Diamond$ and the Cantor derivative operator on topological scattered spaces, and showed that $\gl$ is strongly complete with respect to the class of all scattered spaces.

The study of the completeness of $\gl$ with respect to topologies on ordinals was initiated by Abashidze \cite{abashidze1985ordinal} and Blass \cite{blass}.
Blass interpreted modal operators via filters associated with uncountable cardinals.
He proved that, within ZFC, $\gl$ is complete with respect to end-segment filters. This result implies the completeness of $\gl$ on every ordinal $\Theta\geq \omega^\omega$ equipped with the interval topology, a result independently proved by Abashidz.
Furthermore, Blass proved the completeness of $\gl$ with respect to club filters under suitable square principles $\Box_\kappa$, for all uncountable cardinals $\kappa<\aleph_\omega$.
He also showed that the failure of completeness for club filters is equiconsistent with
the existence of a Mahlo cardinal.
Using the same labelling methods, the completeness of $\gl$ for the filter sequence of normal measures was also established in \cite{GZ}.

The strong completeness of $\gl$ with respect to a specific topological space was first investigated in \cite{aguilera2017strong}. The authors showed that $\gl$ is strongly complete with respect to $(\Theta,\mathcal{I}_\lambda)$ for each ordinal $\Theta>e^\lambda\omega$, where $\mathcal{I}_\lambda$ denotes the generalised Icard topology. In particular, $\gl$ is strongly complete with respect to every ordinal $\Theta>\omega^\omega$ with its interval topology $\mathcal{I}_1$.
The methods used in \cite{aguilera2017strong} rely on the fact that the set of propositional variables $\mathbb{P}$ is countable, and they raised the question of whether $\gl(\mathbb{P})$ is strongly complete with respect to some natural topological space whenever $|\mathbb{P}|\geq\aleph_1$.

The strong completeness of $\gl(\mathbb{P})$ for uncountable languages was further investigated in \cite{GZS}. For any ordinal $\lambda>0$, by using the Erd\H{o}s--Rado theorem, the authors showed the failure of strong completeness of $\gl(\mathbb{P})$ for languages of cardinality $(2^{|\lambda|+\aleph_0})^+$ with respect to ordinal spaces equipped with the generalised Icard topologies $\mathcal{I}_\lambda$.
They also introduced a class of topological spaces, called $\kappa$-bouquet spaces, and proved that $\gl(\mathbb{P})$ is strongly complete with respect to these spaces whenever $|\mathbb{P}|=\kappa\geq\aleph_0$.
However, it remains open whether there is a single space with respect to which $\gl(\mathbb{P})$ is strongly complete.

In this paper, we continue investigating the strong completeness of $\gl(\mathbb{P})$ for uncountable languages. More specifically, we study the strong completeness of $\gl(\mathbb{P})$ when $\mathbb{P}=\{p_\xi : \xi<\kappa\}$ has uncountable size $\aleph_0<\kappa\leq 2^{\aleph_0}$.
We show that for any infinite cardinal $\kappa$, the set
\[
\Gamma_\kappa=\set{\Diamond p_\xi:\xi<\kappa}\cup \set{\Box(p_\xi\to\Diamond p_\eta):\xi<\eta<\kappa}
\]
has a model on an ordinal space with the interval topology if and only if $\kappa\nrightarrow_{\rm ipath}(\omega)^2_{\omega,<\omega}$. The latter notation means that there is a colouring $c:[\kappa]^2\to \omega$ such that the consecutive edges of every infinite increasing path receive infinitely many distinct colours.
We then consider $\gl(\mathbb{P})$ when $|\mathbb{P}|<\mathfrak{p}$, where $\mathfrak{p}$ denotes the pseudo-intersection number, and prove that it is strongly complete with respect to $(\omega^\omega+1,\mathcal{I}_1)$.
Finally, using results from first-order model theory, we prove the strong completeness of $\gl(\mathbb{P})$ for $|\mathbb{P}|=\aleph_1$.

\section{Preliminaries}

In this section, we briefly review some basic notions and facts from set theory and the theory of ordinal models. Most of the material is taken from \cite{GZS}. For further details on topics such as the Kripke semantics of $\gl$ and its ordinal spaces, we refer the reader to \cite{aguilera2017strong} as well.

For an ordinal $\Theta$, let $\mathcal I_1$ denote the interval topology generated by all sets of the form $(\alpha,\beta]$ and $[0,\beta]$.
Over this topology, any limit ordinal is a limit point.
Clearly, $(\Theta,\mathcal{I}_1)$ is a scattered space.
For any $A\subseteq\Theta$, denote the set of its limit points by $d_{\mathcal{I}_1}(A)$.

A model $\mathfrak{M}$ based on $(\Theta,\mathcal{I}_1)$ is obtained by adding a valuation function $\lbr\cdot\rbr$ assigning to each propositional variable $p\in\mathbb{P}$ a subset of ordinals in $\Theta$.
This valuation extends to all formulas by respecting all Boolean connectives and interpreting the modal operator $ \Diamond$ as the derived set operator $d_{\mathcal{I}_1}$.

For a model $\model$ based on $(\Theta,\mathcal I_1)$ and $A\subseteq\Theta$, we write
$
\model,A\models\varphi
$
if $\model,\alpha\models\varphi$ for all $\alpha\in A$.

\begin{fact}[\cite{GZS}]\label{fact:interval}
Let $\model$ be a model based on $(\Theta,\mathcal I_1)$. For any $\alpha<\Theta$ and formula $\varphi$:
\begin{enumerate}[(1)]
\item $\model,\alpha\models\Diamond\varphi$ if and only if
$\set{\beta<\alpha:\model,\beta\models\varphi}$
is cofinal in $\alpha$.
\item If $\tup{\alpha_i:i<\cof(\alpha)}$ is increasing and cofinal in $\alpha$, then
$\model,\alpha\models\Box\varphi$
if and only if there is $i<\cof(\alpha)$ such that
$\model,(\alpha_i,\alpha)\models\varphi$.
\end{enumerate}
\end{fact}

For every non-zero ordinal $\xi$, let $\ell(\xi)$ denote its end logarithm; that is, $\ell(\xi)=\rho$ if and only if $\xi=\gamma+\omega^\rho$ for some ordinal $\gamma$. Also, set $\ell(0)=0$. It is evident that if $\delta>0$, then
$\ell(\gamma+\delta)=\ell(\delta).$
For any set $A$ and cardinal $\kappa$, let $[A]^\kappa$ denote the set of all subsets of $A$ of cardinality $\kappa$.

In \cite{GZS}, by defining a suitable colouring function and applying the Erd\H{o}s--Rado theorem, it is shown that $\Gamma_\kappa$ has no ordinal model whenever $\kappa= (2^{\aleph_0})^+$. In this paper,
we use the following path partition relation to obtain an equivalent condition for the satisfiability of $\Gamma_\kappa$ for any infinite cardinal $\kappa$.
This relation expresses the property that every colouring admits an infinite increasing path whose consecutive edges receive only finitely many distinct colours.

\begin{dfn}\label{def:ipath}
	Let $\kappa$ and $\theta$ be infinite cardinals. We write
	$\kappa\rightarrow_{\rm ipath}(\omega)^2_{\theta,<\omega}$
	to denote the assertion that for every colouring
	$c:[\kappa]^2\longrightarrow\theta$,
	there is a strictly increasing sequence $\langle\xi_n:n<\omega\rangle$ such that
	$\{c(\{\xi_n,\xi_{n+1}\}):n<\omega\}$
	is finite.
\end{dfn}
To ease notation, throughout the paper, we denote $c(\{\xi,\eta\})$ simply by $c(\xi,\eta)$.

Following \cite{Scheepers}, let $M(\theta, 2)$ denote the least cardinal $\kappa$ such that $\kappa\rightarrow_{\rm ipath}(\omega)^2_{\theta,<\omega}$.
Then we have the following; see \cite{Scheepers}, where clause (1) is proved, and clause (2) is attributed to an unpublished result of
Todor\v{c}evi\'c.
\begin{lem}
\label{lem:path}
\begin{enumerate}
\item $M(\aleph_0, 2) > \aleph_1$,
\item It is consistent that $M(\aleph_0, 2)=\aleph_2 < \mathfrak c= 2^{\aleph_0}$.
\end{enumerate}
\end{lem}

Let $(P,\leq)$ be a partially ordered set.
A subset $A\subseteq P$ is \textit{centred} if for each finite subset $F\subseteq A$ there is $p\in P$ such that $p\leq q$ for all $q\in F$.
A poset $(P,\leq)$ is \textit{$\sigma$-centred} if it is the union of countably many centred subsets.

\begin{dfn}\label{def:p}
\begin{enumerate}
\item A family $\mathcal A\subseteq[\omega]^\omega$ has the \textit{strong finite intersection} property if
$
\left|\bigcap\mathcal F\right|=\aleph_0
$
for every $\mathcal F\in[\mathcal A]^{<\omega}$.
\item An infinite $X\subseteq\omega$ is a \textit{pseudo-intersection} of $\mathcal A$ if
 $X\substar A$ for every $A\in\mathcal A$, that is, if $X\setminus A$ is finite for every $A\in\mathcal A$.
\item The cardinal $\pp$, called the \textit{pseudo-intersection number}, is the least cardinality of a family of infinite subsets of $\omega$ that has the strong finite intersection property but has no pseudo-intersection.
\end{enumerate}
\end{dfn}
Consequently, $\mathfrak p$ is the smallest size of a family of infinite subsets of $\omega$ such that every finite subfamily has an infinite intersection, while no infinite set is almost contained in every member of the family.
Note that we have $\aleph_1\leq\pp\leq \mathfrak{c}$. Clearly, if CH holds, $\pp=\aleph_1= \mathfrak{c}$. However, without CH, it is consistent that
$\aleph_1<\mathfrak p<\mathfrak c$.

We use the following form of Bell's theorem for characterising the pseudo-intersection number $\pp$.

\begin{thm}[Bell \cite{Bell81}]\label{thm:bell}
Let $\kappa$ be an infinite cardinal. The following are equivalent.
\begin{enumerate}[(1)]
\item $\kappa<\pp$.
\item If $\mathbb Q$ is a $\sigma$-centred partial order and $\mathcal D$ is a family of at most $\kappa$ dense subsets of $\mathbb Q$, then there is a filter $G\subseteq\mathbb Q$ such that
$
G\cap D\ne\varnothing
$
for every $D\in\mathcal D$.
\end{enumerate}
\end{thm}

\section{Ordinal Models of \(\Gamma_\kappa\)}

Let $\kappa$ be an infinite cardinal, and
$
\mathbb P=\set{p_\xi:\xi<\kappa}
$.
Set
\[
\Gamma_\kappa=
\set{\Diamond p_\xi:\xi<\kappa}
\cup
\set{\Box(p_\xi\to\Diamond p_\eta):\xi<\eta<\kappa}.
\]
Clearly, $\Gamma_\kappa$ is $\gl$-consistent.
In \cite{GZS} it is shown that if $\kappa= \mathfrak{c}^+$, then $\Gamma_\kappa$ has no model of the form $(\Theta,\mathcal{I}_1)$ for any ordinal $\Theta$.
In this section we investigate whether $\Gamma_\kappa$ can have an ordinal model.

\begin{lem}\label{lem:rank-layers}
Assume $0<\theta$, $\rho<\theta$, and set
$
L^\theta_\rho=\set{\delta:0<\delta<\omega^\theta\text{ and }\ell(\delta)=\rho}.
$
Then:
\begin{enumerate}[(1)]
\item $L^\theta_\rho$ is cofinal in $\omega^\theta$.
\item If $\sigma<\rho<\theta$, then
$L^\theta_\rho\subseteq d_{\mathcal I_1}(L^\theta_\sigma)$.
\item For every ordinal $\gamma$, the set
$
\gamma+L^\theta_\rho=\set{\gamma+\delta:\delta\in L^\theta_\rho}
$
is cofinal in $[\gamma,\gamma+\omega^\theta)$. Also if $\sigma<\rho<\theta$, then
$\gamma+L^\theta_\rho
\subseteq
 d_{\mathcal I_1}(\gamma+L^\theta_\sigma)$.
\end{enumerate}
\end{lem}
\begin{proof}
	For (1), fix $\rho<\theta$ and let $\alpha =\sup L^\theta_\rho$.
	If $\alpha<\omega^\theta$, then, since $\omega^\theta$ is additively indecomposable, we have $\beta=\alpha+\omega^\rho<\omega^\theta$. Also, $\ell(\alpha+\omega^\rho)=\rho$. So, $\beta\in L^\theta_\rho$, a contradiction with $\alpha$ being a supremum.
	
	For (2), assume that $\delta\in L^\theta_\rho$ and $(\alpha,\delta]\in\mathcal{I}_1$.
	Since $\alpha<\delta$, we clearly have $\gamma=\alpha+\omega^\sigma<\delta$ and $\gamma\in L^\theta_\sigma$.

	Clause (3) follows from (1) and (2).
\end{proof}

\begin{thm}\label{thm:path-characterization}
Let $\kappa$ be an infinite cardinal. The following are equivalent.
\begin{enumerate}[(1)]
\item $\Gamma_\kappa$ has a model based on an ordinal with the interval topology.
\item $\kappa\nrightarrow_{\rm ipath}(\omega)^2_{\omega,<\omega}$.
\end{enumerate}
Furthermore, if (2) holds, then $\Gamma_\kappa$ has a model based on an ordinal $\Theta <\omega^{\kappa^+}$.
\end{thm}

\begin{proof}
Assume (1). Let $\model$ be a model based on $(\Theta,\mathcal I_1)$, and suppose that
$
\model,\alpha\models\Gamma_\kappa.
$
Then, by \cite[Proposition~3.5]{GZS}, $\alpha$ is a limit ordinal, and it has countable cofinality.
Fix an increasing cofinal sequence
$
\tup{\alpha_n:n<\omega}
$
in $\alpha$. For $\xi<\eta<\kappa$, let
\[
c(\xi,\eta)=\min\set{n<\omega:\model,(\alpha_n,\alpha)\models p_\xi\to\Diamond p_\eta}.
\]
We show that $c$ witnesses
$
\kappa\nrightarrow_{\rm ipath}(\omega)^2_{\omega,<\omega}.
$
Suppose, towards a contradiction, that there is an increasing sequence
$
\xi_0<\xi_1<\cdots
$
such that
$
F=\set{c(\xi_n,\xi_{n+1}):n<\omega}
$
is finite. Let $N=\max F$. Then
\[
\model,(\alpha_N,\alpha)\models p_{\xi_n}\to\Diamond p_{\xi_{n+1}}
\]
for all $n<\omega$.  Using this,  we obtain a strictly decreasing sequence of ordinals, thereby reaching a contradiction (see \cite[Lemma~3.3]{GZS}).

Assume now (2) holds, and let
$
c:[\kappa]^2\longrightarrow\omega
$
witness it. For each $n<\omega$, define a relation $R_n$ on $\kappa$ by
\[
\xi R_n\eta\quad\Longleftrightarrow\quad
\xi<\eta\text{ and }c(\xi,\eta)\le n.
\]
There is no infinite $R_n$-chain; hence $R_n$ is converse well-founded. Define recursively
\[
r_n(\xi)=\sup\set{r_n(\eta)+1:\xi R_n\eta}
\]
and set
$
\theta_n=\sup_{\xi<\kappa}(r_n(\xi)+1).
$
Then
\[
\xi R_n\eta\quad\Longrightarrow\quad r_n(\xi)>r_n(\eta).
\]
If $r_n(\xi)\ge\kappa^+$ for some $\xi<\kappa$, then, by the definition of $r_n$,
there is $\eta$ with $\xi R_n\eta$ and $r_n(\eta)\ge\kappa^+$. Repeating this gives an infinite $R_n$-chain, which is a contradiction. Hence
$
r_n(\xi)<\kappa^+
$
for all $\xi<\kappa$. It follows that
$
\theta_n<\kappa^+,$ and hence $
\theta=\sup_{n<\omega}\theta_n<\kappa^+.
$

Now define
$
\gamma_0=0,
$
$\gamma_{n+1}=\gamma_n+\omega^{\theta_n},
$
and set
$\alpha=\sup_{n<\omega}\gamma_n.$ Clearly
$
\alpha\le\omega^\theta\cdot\omega<\kappa^+.
$
For $\xi<\kappa$, set
$
\inter{p_\xi}=
\bigcup_{n<\omega}S_{n,\xi},
$
where
\[
S_{n,\xi}=\set{\gamma_n+\delta:
0<\delta<\omega^{\theta_n}
\text{ and }
\ell(\delta)=r_n(\xi)}.
\]
Let $\model$ be the corresponding model on $(\alpha+1,\mathcal I_1)$.

By Lemma~\ref{lem:rank-layers}, $S_{n,\xi}$ is cofinal in
$
[\gamma_n,\gamma_{n+1})
$
for every $n<\omega$. Hence $\inter{p_\xi}$ is cofinal in $\alpha$, and therefore
$
\model,\alpha\models\Diamond p_\xi
$
for all $\xi<\kappa$.

Fix $\xi<\eta<\kappa$, and let
$
N=c(\xi,\eta).
$
If $n\ge N$, then $\xi R_n\eta$, and hence
$
r_n(\xi)>r_n(\eta).
$
By Lemma~\ref{lem:rank-layers}(3),
$
S_{n,\xi}\subseteq d_{\mathcal I_1}(S_{n,\eta}).
$
It follows that
$
\model,(\gamma_N,\alpha)\models p_\xi\to\Diamond p_\eta.
$
By Fact~\ref{fact:interval},
$
\model,\alpha\models\Box(p_\xi\to\Diamond p_\eta).
$
Thus
$
\model,\alpha\models\Gamma_\kappa.
$
\end{proof}
By Lemma \ref{lem:path} and Theorem~\ref{thm:path-characterization}, we have the following conclusion.
\begin{cor}\label{cor:omega1-Gamma}
\begin{enumerate}
\item $\Gamma_{\omega_1}$ has a model based on an ordinal with the interval topology.

\item It is consistent that $\mathfrak c > \aleph_2$ and  $\Gamma_{\omega_2}$ has no ordinal model  with the interval topology.

\end{enumerate}
\end{cor}

\section{Strong completeness below \texorpdfstring{$\mathfrak p$}{p}}

Let $\mathbb P$ be a set of propositional variables with
 $
|\mathbb P|\le\mc$
and fix an injection
 $e:\mathbb P\longrightarrow 2^\omega$.
Throughout this paper, we assume that all trees are transitive and converse well-founded.
For a finite rooted tree $T$, $n<\omega$, and a map
 $
a:2^n\times T\longrightarrow2,
 $
define a valuation $V_{n,a}$ on $T$ by
\[
V_{n,a}(p) = \{t\in T\; : \; a(e(p)\restr n,t)=1 \}.
\]

Let $\mathscr D$ be the collection of all pointed models of the form
 $
(T,V_{n,a},\root(T)),
 $
where $T$ ranges over a fixed countable set of representatives of the finite rooted trees.
For $\varphi\in\lang(\mathbb P)$, set
\[
A_\varphi=\set{d\in \mathscr D:d\models\varphi}.
\]
where by $d\models\varphi$ we mean that $d,\root(T)\models\varphi$ when $d=(T,V_{n,a},\root(T))$.

\begin{lem}\label{lem:universal-family}
Let $\Sigma\subseteq\lang(\mathbb P)$ be finite and $\gl$-consistent. Then
$
\bigcap_{\varphi\in\Sigma}A_\varphi
$
is infinite.
\end{lem}

\begin{proof}
	Let $M=(T,V,r)$
	be  a finite rooted model of
	$\Sigma.$
	Let
	$
	\set{p_0,\ldots,p_{m-1}}
	$
	be the set of variables occurring in $\Sigma$. Choose $N<\omega$ such that
	$
	e(p_i)\restr N\ne e(p_j)\restr N
	$
	for all $i<j<m$. 
	Given $n \geq N$, let $a_n:2^n\times T\to2$ be such that
	\[
	a_n(e(p_i)\restr n,t)=1
	\quad\Longleftrightarrow\quad
	t\in V(p_i)
	\]
	for all $i<m$ and $t\in T$. Then
	$
	(T,V_{n,a_n}),r\models\Sigma
	$.
\end{proof}

Let $\Gamma\subseteq\lang(\mathbb P)$ be maximal $\gl$-consistent, and set
\[
\mathcal B_\Gamma=
\set{\varphi\wedge\Box\varphi:\Box\varphi\in\Gamma}
\]
and
\[
\mathcal E_\Gamma=
\set{\psi:\Diamond\psi\in\Gamma}.
\]

\begin{lem}\label{lem:finite-compatibility}
If $F\in[\mathcal B_\Gamma]^{<\omega}$ and $\psi\in\mathcal E_\Gamma$, then
$
F\cup\set{\psi}
$
is $\gl$-consistent.
\end{lem}
\begin{proof}
The claim holds since $\Gamma$ is consistent and $\Gamma\vdash \Diamond(\bigwedge F\wedge\psi)$.
\end{proof}

We now define a forcing notion $\mathbb Q_\Gamma$ on $[\mathscr D]^{<\omega}\times[\mathcal B_\Gamma]^{<\omega}$ as follows:
 \begin{itemize}
\item  a condition is a pair
$
p=(s_p,F_p)
$
with
$
s_p\in[\mathscr D]^{<\omega}$ and
$F_p\in[\mathcal B_\Gamma]^{<\omega}$;

\item for $p,q\in\mathbb Q_\Gamma$, let  $q\le p$ ($q$ is stronger than $p$) if
$
s_p\subseteq s_q,$
$
F_p\subseteq F_q,
$
and
$
s_q\setminus s_p\subseteq\bigcap_{b\in F_p}A_b.
$
\end{itemize}
\begin{lem}\label{lem:sigma-centered}
$\mathbb Q_\Gamma$ is $\sigma$-centred.
\end{lem}

\begin{proof}
For $s\in[\mathscr D]^{<\omega}$ set
\[
C_s=\set{(s,F):F\in[\mathcal B_\Gamma]^{<\omega}}.
\]
Each $C_s$ is centred.
Hence, as $\mathscr{D}$ is countable,
$\mathbb Q_\Gamma=\bigcup_{s\in[\mathscr{D}]^{<\omega}}C_s$
is $\sigma$-centred.
\end{proof}

For $b\in\mathcal B_\Gamma$, $\psi\in\mathcal E_\Gamma$, and $m<\omega$, set
\[
\mathcal D_b=\set{p\in\mathbb Q_\Gamma:b\in F_p}
\]
and
\[
\mathcal E_{\psi,m}=
\set{p\in\mathbb Q_\Gamma:|s_p\cap A_\psi|\ge m}.
\]

\begin{lem}\label{lem:dense-sets}
The sets $\mathcal D_b$ and  $\mathcal E_{\psi,m}$ are dense in $\mathbb Q_\Gamma$.
\end{lem}

\begin{proof}
Density of $\mathcal D_b$ is immediate. Let $p=(s,F)$ and let $\psi\in\mathcal E_\Gamma$. By Lemmas~\ref{lem:finite-compatibility} and \ref{lem:universal-family},
$
A_\psi\cap\bigcap_{b\in F}A_b
$
is infinite. Let $t$ be a set consisting of $m$ elements of the above set and put $q=(s \cup t, F)$. Then $q$ is an extension of $p$ in
$\mathcal E_{\psi,m}$.
\end{proof}
The following proposition is an immediate consequence of the above lemma.
\begin{prop}\label{prop:selected-X}
Assume $G\subseteq\mathbb Q_\Gamma$ is a filter meeting every $\mathcal D_b$ and every $\mathcal E_{\psi,m}$. Set
$
X_G=\bigcup_{p\in G}s_p.
$
Then
$
X_G\substar A_b
$
for all $b\in\mathcal B_\Gamma$, and
$
|X_G\cap A_\psi|=\aleph_0
$
for all $\psi\in\mathcal E_\Gamma$.
\end{prop}

We use the following consequence of the ordinal $d$-map construction; see \cite[Lemma~6.1]{aguilera2017strong}.

\begin{fact}\label{fact:finite-block}
Let $T$ be a finite rooted tree with the root $t$. There are an ordinal $\delta<\omega^\omega$ and a surjective $d$-map
\[
f:(\delta+1,\mathcal I_1)\longrightarrow T
\]
such that
$
f(\delta)=t.
$
\end{fact}

\begin{lem}\label{lem:selection-to-ordinal}
Let $\Gamma$ be maximal $\gl$-consistent and suppose that there is an infinite $X\subseteq \mathscr{D}$ such that
$
X\substar A_b
$
for every $b\in\mathcal B_\Gamma$, and
$
|X\cap A_\psi|=\aleph_0
$
for every $\psi\in\mathcal E_\Gamma$. Then $\Gamma$ is satisfiable on
$
(\omega^\omega+1,\mathcal I_1).
$
\end{lem}

\begin{proof}
Let
$
\set{\mathfrak{M}_n:n<\omega}
$
be an enumeration of $X$,
 where
$
\mathfrak{M}_n=(T_n,V_n,t_n).
$
By Fact~\ref{fact:finite-block}, for each $n<\omega$, choose $\delta_n<\omega^\omega$ and a surjective $d$-map
\[
f_n:(\delta_n+1,\mathcal I_1)\longrightarrow T_n
\]
with $f_n(\delta_n)=t_n$.
Define
\[
\beta_0=0,
\qquad
\beta_{n+1}=\beta_n+\delta_n+1,
\qquad
\alpha=\sup_{n<\omega}\beta_n.
\]
Also, set
$
B_n=[\beta_n,\beta_n+\delta_n].
$
Each $B_n$ is open in $(\alpha+1,\mathcal I_1)$ and is order-homeomorphic to $\delta_n+1$. Denote this homeomorphism by $i_n$.
Furthermore, by additive indecomposability of $\omega^\omega$, it is easy to see $\alpha\le\omega^\omega.$

Now, define a valuation $\lbr\cdot\rbr$ over $\alpha+1$ as the union of all $(f_n\circ i_n)^{-1}(V_n)$ plus putting, at $\alpha$,
\[
\alpha\in\inter p
\quad\Longleftrightarrow\quad
p\in\Gamma.
\]
Let $\mathfrak{M}=(\alpha+1,\mathcal{I}_1,\lbr\cdot\rbr)$.
We prove by induction on $\varphi$ that
\[
\model,\alpha\models\varphi
\quad\Longleftrightarrow\quad
\varphi\in\Gamma.
\]
The atomic and Boolean cases are immediate.
Suppose $\Box\varphi\in\Gamma$. Set
$
b=\varphi\wedge\Box\varphi \in \mathcal B_\Gamma.
$
Then, by the assumption $
X\substar A_b
$, we have $\mathfrak{M}_n,t_n\models b$ for all but finitely many $n$. Since
$
\mathfrak{M}_n,t_n\models\varphi\wedge\Box\varphi$,
every point of $T_n$ satisfies $\varphi$. Hence every point of $B_n$ satisfies $\varphi$. Therefore there is $N<\omega$ such that
$
\model,(\beta_N,\alpha)\models\varphi.
$
By Fact~\ref{fact:interval},
$
\model,\alpha\models\Box\varphi.
$

Suppose now $\Box\varphi\notin\Gamma$. Then
$
\Diamond\neg\varphi\in\Gamma,
$
hence $\neg\varphi \in \mathcal E_\Gamma$.
By the assumption $
|X\cap A_{\neg\varphi}|=\aleph_0
$,
there are infinitely many $n<\omega$ such that
$
\mathfrak{M}_n,t_n\models\neg\varphi.
$
For each such $n$,
$
\model,\beta_n+\delta_n\models\neg\varphi.
$
These points are cofinal in $\alpha$, hence
$
\model,\alpha\models\Diamond\neg\varphi.
$

Thus $\Gamma$ is satisfiable at $\alpha$. Since $[0,\alpha]$ is open in $(\omega^\omega+1,\mathcal I_1)$, extend the valuation arbitrarily to the remaining points of $\omega^\omega+1$.
\end{proof}

\begin{thm}\label{thm:below-p}
Suppose $\kappa<\pp$ is infinite and
$
|\mathbb P|\le\kappa.
$
Then $\gl(\mathbb P)$ is strongly complete with respect to
$
(\omega^\omega+1,\mathcal I_1).
$
\end{thm}

\begin{proof}
Let $\Gamma\subseteq\lang(\mathbb P)$ be $\gl$-consistent; without loss of generality we can assume that it is maximal $\gl$-consistent. If $\Box\bot\in\Gamma$, then it is satisfiable on an isolated point. Suppose $\Diamond\top\in\Gamma$.
 By Lemma~\ref{lem:sigma-centered}, $\mathbb Q_\Gamma$ is $\sigma$-centred. The family
\[
\set{\mathcal D_b:b\in\mathcal B_\Gamma}
\cup
\set{\mathcal E_{\psi,m}:\psi\in\mathcal E_\Gamma,\ m<\omega}
\]
has cardinality at most $\kappa <\pp$, hence by Theorem~\ref{thm:bell}, there is a filter $G\subseteq\mathbb Q_\Gamma$ meeting every member of this family. By Proposition~\ref{prop:selected-X}, $X_G$ satisfies the assumptions of Lemma~\ref{lem:selection-to-ordinal}. Hence $\Gamma$ is satisfiable on $(\omega^\omega+1,\mathcal I_1)$.
\end{proof}

Note that since $\aleph_0 < \pp$, the above theorem provides another proof of the result of \cite{aguilera2017strong}.

\begin{cor}\label{cor:path-below-p}
If $\kappa<\pp$ is infinite, then
$
\kappa\nrightarrow_{\rm ipath}(\omega)^2_{\omega,<\omega}.
$
In particular $M(\aleph_0, 2) \geq \pp.$
\end{cor}

Recall that Martin's axiom, denoted by $\mathsf{MA}$, says that if $P$ is a ccc poset, $\kappa < \mathfrak{c}$, and
$\mathcal{D}=\{D_\alpha : \alpha<\kappa\}$
is a family of $\kappa$ dense subsets of $P$, then there is a filter $G\subseteq P$ such that
$G\cap D_\alpha\neq \emptyset$
for every $\alpha<\kappa$.
Let
$\mathsf{MA}(\sigma\text{-centred})$ denote the version of this axiom in which the poset is required to be $\sigma$-centred rather than merely ccc.
The following corollary is obtained by using Theorem~\ref{thm:bell} and Corollary~\ref{cor:path-below-p}.
\begin{cor}\label{cor:ma-omega1}
Assume $\mathsf{MA}(\sigma\text{-centred})+\mathfrak c > \aleph_1$. If $|\mathbb P| < \mathfrak c$, then $\gl(\mathbb P)$ is strongly complete with respect to
$
(\omega^\omega+1,\mathcal I_1).
$
\end{cor}

\section{Strong Completeness for Languages of Size $\aleph_1$}

In this section, by using some model-theoretic results, we study the strong completeness of $\gl(\mathbb{P})$ for $|\mathbb{P}|=\aleph_1$; for more details on model theory, see e.g., \cite{marker}.

Recall that $\gl$ is complete with respect to the class of finite trees.
Any finite rooted tree $(T,<_T)$ can be seen as a modal algebra
$\mathcal{A}_T=(\mathcal{P}(T), 0, 1, \wedge,\neg ,d)$, where
\[
d(X)=\{w\in T : \exists v\in X \; \text{s.t.}\; wRv\},
\]
for each $X\subseteq T$.
Note that these algebras are a specific kind of Magari algebras that correspond to the scattered spaces (see, e.g., \cite{beklemishev2014}).

Let
\[\mathcal{K}=\{\mathcal{A}_T \ : \ (T,<_T) \text{ is a finite rooted tree}\}.
\]
Any $\mathcal{A}_T\in\mathcal{K}$ can be seen as a first-order structure in the language
$\mathcal{L}^\star=\{0,1,\wedge,\neg,d\}$. Furthermore, we can assign to each modal formula $\varphi$ a first-order formula $\varphi^\star$ in the language $\mathcal{L}^\star_1=\mathcal{L}^\star\cup\{c_\xi : \xi<\omega_1\}$,  where $c_\xi$ is a constant symbol (and is obtained by replacing each occurrence of $p_\xi$ with $c_\xi$ when translating $\varphi$ into $\varphi^\star$).
 In this case, any model $\mathfrak{M}=(T,<_T,V)$ can be seen as an extension of $\mathcal{A}_T$ in the language $\mathcal{L}^\star_1$, and vice versa.
	Thus, a modal formula $\varphi$ is satisfiable in $(T,<_T)$ if and only if $\varphi^\star\neq 0$ is satisfiable in some $\mathcal{L}^\star_1$-extension of $\mathcal{A}_T$. 

Suppose that $|\mathbb{P}|=\aleph_1$ and let $\Gamma$ be a maximal $\gl(\mathbb{P})$-consistent set of formulas.
Let
\[
\Sigma_\Gamma =\{\varphi^\star =1 : \Box\varphi \in\Gamma\} \cup \{\psi^\star \neq 0 :  \Diamond\psi \in\Gamma\}.
\]
If $ \Diamond\top\in\Gamma$, then $\Sigma_\Gamma$ is consistent; since in this case the set $\{\varphi\wedge\Box\varphi: \Box\varphi\in \Gamma\} \cup \{ \Diamond\psi : \Diamond\psi\in\Gamma\}$ is consistent.

The following proposition plays a key role in the proof of the completeness theorem.
\begin{prop}\label{prop::reduced-product}
	Let $\Gamma$ be a maximal consistent set of modal formulas and $ \Diamond\top\in\Gamma$. Then there are countably many finite rooted trees  $(T_n,<_n)\in\mathcal{K}$  such that $\Sigma_\Gamma$ is satisfiable in some $\mathcal{L}^\star_1$-extension of $\prod_{n<\omega}\mathcal{A}_{T_n}/\Fin$,
	where $\prod_{n<\omega}\mathcal{A}_{T_n}/\Fin$ is the reduced product modulo the cofinite filter $\Fin$.
\end{prop}
Before proving this proposition, we recall some model-theoretic results from first-order logic.
First, observe that the following equivalence holds for every first-order $\mathcal{L}$-atomic formula $\delta$ and any countable family of $\mathcal{L}$-structures $\{\mathcal{M}_n : n<\omega\}$:
\[
\{n<\omega : \mathcal{M}_n\models\delta\}\text{ is cofinite}
\quad\Longleftrightarrow\quad
\prod_{n<\omega} \mathcal{M}_n/\Fin\models\delta.
\]
Furthermore, it is known (see \cite{JO68}) that
	if $\mathcal{L}$ is a countable first-order language  and  $\{\mathcal{M}_n: n<\omega\}$ is a countable family of $\mathcal{L}$-structures, then the reduced product $\prod_{n<\omega}\mathcal{M}_n/\Fin$ is $\aleph_1$-saturated. 

For an $\mathcal{L}$-structure $\mathcal{M}$, let  $\diag(\mathcal{M})$ be the atomic diagram of
$\mathcal{M}$,  in the language
$\mathcal{L}_{\mathcal{M}}$.

\begin{lem}\cite[Lemmas 2.3.3, and ~4.3.17]{marker}\label{lem:marker}
	\begin{enumerate}
		\item Assume $\mathcal{N}$ is an $\mathcal{L}_{\mathcal{M}}$-structure and it is a model of $\diag(\mathcal{M})$.
		Then, viewing $\mathcal{N}$ as an $\mathcal{L}$-structure, there is an embedding $j: \mathcal{M} \to \mathcal{N}$.
		
		\item 	Let $\kappa\geq\aleph_0$. If $\mathcal{M}$ is $\kappa$-saturated, then it is $\kappa^+$-universal in the sense that if $\mathcal{N}\equiv\mathcal{M}$ and $|N|\leq\kappa$, then there is an elementary embedding $e:N\to M$.
	\end{enumerate}
\end{lem}

By an \textit{algebraic language}, we mean a first-order language $\mathcal{L}$ that contains no relational symbols. We also call any structure of this language an \textit{$\mathcal{L}$-algebra}.
\begin{lem}\label{lem:embedding-reduced-product}
	Assume $\mathcal{L}$ is  a countable algebraic language,  $\mathcal{A}$ is an
	$\mathcal{L}$-algebra of size $|A|\leq\aleph_1$, and let
	$\langle \mathcal{B}_n:n<\omega \rangle$ be a sequence of $\mathcal{L}$-algebras.
	Assume for every finite conjunction $\delta(\bar{x})$ of equations
	and inequations, if
	$\mathcal{A}\models\exists\bar{x}\,\delta(\bar{x}),$
	then for all sufficiently large $n<\omega,$
	$\mathcal{B}_n\models\exists\bar{x}\,\delta(\bar{x})$.
	Then there is an $\mathcal{L}$-embedding
	$e:\mathcal{A}\to\prod_{n<\omega}\mathcal{B}_n/\Fin.$
\end{lem}
\begin{proof}
	Set
	$\mathcal{R}=\prod_{n<\omega}\mathcal{B}_n/\Fin,$
	and note that it is $\aleph_1$-saturated.
	By our assumption,
	$\Th(\mathcal{R})\cup \diag(\mathcal{A})$
	is consistent.
	 To see this, assume that $\{\delta_1,\cdots,\delta_n\}$ is a finite subset of $\diag(\mathcal{A})$. Then we have $\mathcal{A}\models \exists \bar{x}\delta(\bar{x})$, where $\delta(\bar{x})= \bigwedge_{i\leq n}\delta_i(\bar{x})$, which is obtained by replacing all constants with variables. Then $\mathcal{B}_n\models\exists\bar{x}\delta(\bar{x})$ for all sufficiently large $n$. Hence, since $\delta_i$ is atomic or negated atomic, one can find an appropriate $[\bar{(b_n)}]$ such that $\mathcal{R}\models\delta_i(\bar{[(b_n)]})$. 
	Thus, by compactness and downward L\"owenheim--Skolem, this theory has a model
	$\mathcal{C}$ of cardinality at most $\aleph_1$. Since  $\mathcal{C}\equiv\mathcal{R}$ as $\mathcal{L}$-structures,
	Lemma \ref{lem:marker}(1) gives an $\mathcal{L}$-embedding
	$j:\mathcal{A}\to\mathcal{C}$.
	On the other hand,
	$\mathcal{R}$ is $\aleph_1$-saturated and $|C|\leq\aleph_1$, hence by
	Lemma \ref{lem:marker}(2), there is an elementary embedding
	$k:\mathcal{C}\to\mathcal{R}.$
	Then $e=k\circ j$ is the required $\lang$-embedding.
\end{proof}

\begin{proof}[Proof of Proposition~\ref{prop::reduced-product}]
	First, note that since $\Gamma$ is maximal consistent and $ \Diamond\top\in \Gamma$, $\Sigma_\Gamma$ is consistent in $\mathcal{K}$. 
	
	Let $\Delta =\bigcap\{\Th(\mathcal{A}_T) : \mathcal{A}_T\in\mathcal{K}\}$.
	Then clearly $\Sigma_\Gamma\cup\Delta$ is consistent, so it has a model
	$(\mathcal{A}, \langle c_\xi:\xi<\omega_1\rangle)$ where $\mathcal{A}$ is an $\mathcal{L}^\star$-structure with the domain $A$ such that $|A|\leq\aleph_1$.
	As $\mathcal{L}^\star$ is finite, $\Th_{\mathcal{L}^\star}(\mathcal{A})=\{\theta_n : n<\omega\}$ is countable.
	For each $n<\omega$, let $(T_n,<_n)\in\mathcal{K}$ be such that $\mathcal{A}_{T_n}\models\bigwedge_{i\leq n}\theta_i$.
	Such a $T_n$ exists, since otherwise $\mathcal{A}_T\models\neg\bigwedge_{i\leq n}\theta_i$ for all $(T,<_T)$, so $\neg\bigwedge_{i\leq n}\theta_i\in\Delta$, a contradiction with $\mathcal{A}\models\Delta$.
	
	Thus, by Lemma~\ref{lem:embedding-reduced-product}, there is an $\mathcal{L}^\star$-embedding
	\[
	e:\mathcal{A}\to \prod_{n<\omega}\mathcal{A}_{T_n}/\Fin.
	\]
	Each formula from $\Sigma_\Gamma$ is satisfied in $(\mathcal{A}, \langle c_\xi:\xi<\omega_1\rangle)$ and is an equation or its negation; hence it also holds in the extension $\prod_{n<\omega}\mathcal{A}_{T_n}/\Fin$ obtained by interpreting each $c_\xi$ as $e(c_\xi^\mathcal{A})$.
\end{proof}

\begin{thm}
	Let $|\mathbb{P}|\leq\aleph_1$. Then $\gl(\mathbb{P})$ is strongly complete  with respect to $(\Theta,\mathcal{I}_1)$ for any $\Theta>\omega^\omega$.
\end{thm}
\begin{proof}
	Let $\Gamma$ be a maximal consistent set of $\mathcal{L} (\mathbb{P})$-formulas.
	If $\Box\bot\in\Gamma$, then clearly $\Gamma$ is satisfiable on an isolated point.
	So assume that $ \Diamond\top\in\Gamma$. By Proposition~\ref{prop::reduced-product}, there are finite rooted trees $(T_n,<_n)$, $n<\omega$, such that $\Sigma_\Gamma$ is satisfiable in some extension of $\prod_{n<\omega}\mathcal{A}_{T_n}/\Fin$. 
	So, for each $n<\omega$, one can extract a valuation $V_n$ over $T_n$ from this extension.
	Let $\mathfrak{M}_n=(T_n,<_n,V_n)$.
	Then we have:
	\begin{itemize}
		\item If $\Box\varphi\in\Gamma$, then $\mathfrak{M}_n\models\varphi$ for all large $n$,
		\item if $ \Diamond\varphi\in\Gamma$, then for infinitely many $n$ there exists $w_n\in T_n$ such that $\mathfrak{M}_n,w_n\models\varphi$.
	\end{itemize}
	Now argue as before (see Lemma \ref{lem:selection-to-ordinal}). There are countably many $(T_n,<_n,V_n)$ and for each there is a surjective d-map $f_n:\delta_n+1\to T_n$ (with $\delta_n<\omega^\omega$).
	Thus, by gluing them together as before, we obtain a model of $(\omega^\omega+1,\mathcal{I}_1)$ for $\Gamma$.
\end{proof}

\subsection*{Acknowledgment}
The first author's research has been supported by a grant from IPM (No. 1405030417). This work is also based on research funded by the Iran National Science Foundation (INSF) under project No. 40401358.

\end{document}